\documentclass [12pt]{article}
\usepackage {graphicx, amscd}
\usepackage {longtable, amsthm, amsmath, latexsym, amssymb, epsfig}

\newcommand{\diverg}{\mathop{\rm div}\nolimits}
\newcommand{\Isom}{\mathop{\rm Isom}\nolimits}
\newcommand{\const}{\mathop{\rm const}\nolimits}
\newcommand{\Mor}{\mathop{\rm Mor}\nolimits}

\begin{document}

\begin{center}
\textbf{Factorization of nonlinear heat equation posed on Riemann manifold
\protect\footnote 
{Supported by the grants RFBR 00-15-96042, 99-01-00326.}
}
\end{center}

\begin{center}
Marina F.Prokhorova

Mathematics and Mechanics Institute, Ural Branch RAS, Ekaterinburg, Russia

pmf@imm.uran.ru, http://vpro.convex.ru/Marina

\end{center}

{ 
In [1-5] there was proposed the method of a factorization of PDE. 
The method is based on reduction of complicated systems to more easy ones 
(for example, due to dimension decrease). 
This concept is proposed in general case for the arbitrary 
PDE systems, and its concrete investigation is developing for the heat 
equation case. 

There is considered the category of second order parabolic equations posed 
on arbitrary manifolds. In this category, for the given nonlinear heat 
equation we could find morphisms from it to other parabolic equations with 
the same or a smaller number of independent variables. This allows to receive 
some classes of solutions of original equation from the class of all 
solutions of such a reduced equation. Classification of morphisms (with the 
selection from every equivalence class of the simplest "canonical" 
representatives) is carried out. There are derived the necessary and 
sufficient conditions for canonical morphisms of heat equation to the 
parabolic equation on the other manifold. These conditions are formulated in 
the differential geometry language.

The comparison with invariant solutions classes, obtained by the Lie group 
methods, is carried out. It is proved that discovered solution classes are 
richer than invariant solution classes, even if we find any (including 
discontinuous) symmetry groups of original equation. 
}

\section{General equation category}

\textbf{Definition 1.}
\textit{Task} is a pair $A = \left( {N_{A} ,E_{A}}  \right)$ 
where $N_{A} $ is a set, $E_{A} $ is a system of equations for graph 
$\Gamma \subset N_{A} = M_{A} \times K_{A}$ 
of a function $u: M_{A} \to K_{A} $.

Let $S\left( {A} \right)$ be a set of all subsets $\Gamma \subset N_{A} $ 
satisfying $E_A$. 

\textbf{Definition 2.}
We will say that a (ordered) pair of a tasks $A=\left({N_{A} ,E_{A}} \right)$, 
$B = \left( {N_{B} ,E_{B}}  \right)$ \textit{admits a map} 
$F_{AB} :\;N_{A} \to N_{B} $, if for any 
$\Gamma \subset N_{B}  \quad \Gamma \in S\left( {B} \right) 
\Leftrightarrow F_{AB}^{-1} \left({\Gamma} \right)\in S\left( {A} \right)$. 

Of course, these definitions are rather informal, but they will be correct 
when we define more exactly the notion "system of equations" and the class 
of assumed subsets $\Gamma \subset N_{A} $. Let's consider the 
\textit{general equation category} $\mathcal E$, whose objects are tasks (with some 
refinement of the sense of the notion "system of equations"), and morphisms 
$\Mor \left( {A,B} \right)$ are admitted by the pair $\left( {A,B} \right)$ 
maps with natural composition law. 

For the given task $A$ we could define set $\Mor \left( {A, \mathcal {A}} \right)$ of 
all morphisms $A$ in a framework of some fixed subcategory $\mathcal {A}$ of the 
general equation category (let's call such morphisms and corresponding tasks 
$B$ "factorization of $A$"). The tasks, which factorize $A$, are naturally 
divided into classes of isomorphic tasks, and morphisms $\Mor \left( {A,\; 
\cdot}  \right)$ are divided into equivalence classes. 

The proposed approach is conceptually close to the developed in 
[6] approach to investigation of dynamical and controlled 
systems. In this approach as morphisms of system $A$ to the system $B$ 
there are considered smooth maps 
of the phase space of system $A$ to the phase space of system $B$, which 
transform solutions (phase trajectories) of $A$ to the solutions of $B$. By 
contrast, in the approach presented here, for the class of all solutions of 
reduced system $B$ there is a correspondent class of such solutions of 
original system $A$, which graphs could be projected onto the space 
of dependent and independent variables of $B$; when we pass to the 
reduced system, the number 
of dependent variables remains the same, and the number of independent 
variables does not increase. Thus the approach proposed is an analog to the 
sub-object notion (in terminology of [6]) with respect 
to information about 
original system solutions, though it is closed to the factor-object notion 
with respect to relations between original and reduced systems.

If $G$ is symmetry group of $E_{A} $, then natural projection $p:N \to {{N} 
\mathord{\left/ {\vphantom {{N} {G}}} \right. \kern-\nulldelimiterspace} 
{G}}$ is admitted by the pair $\left( {A,{{A} \mathord{\left/ {\vphantom 
{{A} {G}}} \right. \kern-\nulldelimiterspace} {G}}} \right)$ in the sense of 
Definition 2, that is our definition is the generalization of the reduction 
by the symmetry group. Instead of this general notion of the group analysis 
we base on a more wide notion "a map admitted by the task". We need not 
require from the group preserving solution of an interesting class (if even 
such a group should exist) to be continuous and be admissible by original 
system. So we could obtain more common classes of solutions, than classes of 
invariant solutions of Lie group analysis (though our approach is 
more laborious owing to non-linearity of a system for admissible map). 
Besides, when we factorize original system, defined here factorizing map is 
a more natural object than the group of transformations, operating on space 
of independent and dependent variables of the original task.

\section{Category of parabolic equations} 

Let us consider subcategory $\mathcal {PE}$ of the general equation category, 
whose objects are second type parabolic equations: 
\[
E:u_{t} = Lu,\quad M = T \times X,\quad \quad K = \mathbb R,
\]
\noindent
where $L$ is differential operator, depending on the time $t$, defined on 
the connected manifold $X$, which has the following form in any local 
coordinates $\left( {x^{i}} \right)$ on $X$:
\[ 
Lu = b^{ij}\left( {t,x,u} \right)u_{ij} + c^{ij}\left( {t,x,u} \right)u_{i} 
u_{j} + b^{i}\left( {t,x,u} \right)u_{i} + q\left( {t,x,u} \right).
\]
\noindent
Here a lower index $i$ denotes partial derivative by $x^i$, 
form $b^{ij}= b^{ji}$ is positively defined, $c^{ij}= c^{ji}$. 
Morphisms of $\mathcal {PE}$ are all smooth maps 
admitted by $\mathcal {PE}$ task pairs. Let us describe this morphisms:

\textbf{Theorem 1.} 
Any morphism of the category $\mathcal {PE}$ has the form
\label{Prokhorova:EqMPE}
\begin{equation}
\left( {t,x,u} \right) \to \left( {t'\left( {t} \right),x'\left( {t,x} 
\right),u'\left( {t,x,u} \right)} \right).
\end{equation}
\noindent
Set of isomorphisms of the category $\mathcal {PE}$ is the set of all one-to-one maps 
of kind (1).

Let us consider full subcategory $\mathcal {PE'}$ of the category 
$\mathcal {PE}$, whose objects are equations $u_{t} = Lu$, 
where operator $L$ in local coordinates has the following form:
\[
Lu = b^{ij}\left( {t,x} \right)\left( {a\left( {t,x,u} \right)u_{ij} + 
c\left( {t,x,u} \right)u_{i} u_{j}}  \right) + b^{i}\left( {t,x,u} 
\right)u_{i} + q\left( {t,x,u} \right),
\]
\noindent
and all morphisms are inherited from $\mathcal {PE}$.

\textbf{Theorem 2.} 
If set of morphisms $\Mor _{\mathcal {PE}} \left( {A,B} \right)$ is 
nonempty and $A \in \mathcal {PE'}$, then $B \in \mathcal {PE'}$.

\section{Category of autonomous parabolic equations}

\noindent
Let's call the map (1) \textit{autonomous}, if it has the form 
\label{Prokhorova:EqAutMor}
\begin{equation}
\left( {t,x,u} \right) \to \left( {t,x'\left( {x} \right),u'\left( {x,u} 
\right)} \right).
\end{equation}
\noindent
Let's call a parabolic equation from the category $\mathcal {PE'}$, defined on a 
Riemann manifold $X$, \textit{autonomous}, if it has the form:
\[
u_{t} = Lu = a\left( {x,u} \right)\Delta u + c\left( {x,u} \right)\left( 
{\nabla u} \right)^{2} + \xi \left( {x,u} \right)\nabla u + q\left( {x,u} 
\right),
\quad
\xi \left( { \cdot ,u} \right) \in T^{\ast} X
\]
\textbf{Theorem 3.} 
Let $F:A \to B$ be a morphism of the category $\mathcal {PE}$, $F$ 
be an autonomous map, $A$ be an autonomous equation. Then we could endow 
with Riemann metric the manifold, on which $B$ is posed, in such a way, that 
$B$ becomes an autonomous equation.

Let $\mathcal {APE}$ be the subcategory of $\mathcal {PE'}$, objects of which 
are autonomous parabolic equations, and morphisms are autonomous morphisms 
of the category $\mathcal {PE}$.

\section{Classification of morphisms of nonlinear heat equation}

Let's consider nonlinear heat equation $A \in \mathcal {APE}$, posed on some 
Riemann manifold $X$:
\label{Prokhorova:EqHeat}
\begin{equation}
u_{t} = a\left( {u} \right)\Delta u+q(u).
\end{equation}
\noindent
(note that any equation $u_{t}= a(u)\Delta u+c(u)(\nabla u)^2+q(u))$ 
is isomorphic to some equation (3) in $\mathcal {APE}$). 
We will investigate set of morphisms $\Mor \left( {A,\mathcal {PE}} \right)$ 
and classes of solutions of equation $A$, corresponding these morphisms.

Note, that two morphisms $F:A \to B$ and $F':A \to B'$ are called to be 
equivalent if there exists such isomorphism $G:B \to B'$ that $F' = G \circ 
F$. From the point of view of classes of original task solutions obtained 
from factorization, equivalent morphisms have the same value, that is 
solution classes are the same for these morphisms. So it is interesting to 
select from any equivalence class of the simplest (in some sense) morphism, 
or such morphism for which the factorized equation is the simplest.

When we classify morphisms for the original equation (2), 
a form of coefficient $a\left( {u} \right)$ is important. We will distinguish 
such variants:

\vspace{5mm} 

\quad \quad \quad 
\fbox {$a(u)$ --- arbitrary function}

\vspace{5mm} 

\fbox {$a(u)=a_0 (u - u_0)^{\lambda} $} 
\quad \quad 
\fbox {$a(u)=a_0 e^{\lambda u} $}

\vspace{5mm} 

\quad \quad \quad \quad \quad \quad \quad 
\fbox {$a(u)=a_0$}

\vspace{5mm} 

The lower the variant is situated on this scheme, the richer a collection of 
morphisms is. Note, that similar relation is observed in the group 
classification of nonlinear heat equation [7].

\textbf{Theorem 4.} 
If $a \ne \const$ then for any morphism of equation (3) 
into the category $\mathcal {PE}$ there exists an equivalent in $\mathcal {PE}$
autonomous morphism (that is morphism of the category $\mathcal {APE}$).

Let us give a map $p:X \to X'$ from the manifold $X$ to the manifold $X'$ 
and a differential operator $D$ on $X$. We will say that $D$ \textit{is 
projected} \textit{on} $X'$, if such a differential operator $D'$ on $X'$ 
exists that the following diagram is commutative: 
\[
\begin{CD}
C^{\infty}\left( X'\right) @>{p^{\ast}}>>
C^{\infty}\left( X\right) \\
@V{D'}VV  @VV{D}V \\
C^{\infty}\left( X'\right) @>>{p^{\ast}}>
C^{\infty}\left( X\right) 
\end{CD}
\]

\textbf{Theorem 5.} 
Let $a \ne \const$. For any morphism of the equation $A$ into the 
category $\mathcal {PE}$ there exists an equivalent in $\mathcal {PE}$
autonomous morphism $\left( {t,x,u} \right) \to \left( {t,y\left( {x} 
\right),v\left( {x,u} \right)} \right)$ $A$ to 
$B \in \mathcal {APE}$, for which factorized equation 
$B$ is $v_{t} = a\left( {v} \right)Lv+Q(v)$, 
operator $L$ is projection onto $Y$ at map $x \to y\left( {x} \right)$ of 
the described below operator $D$ (note that this condition is limitation on 
the projection $y\left( {x} \right)$), where:

1) If $A$ is arbitrary (not any of the following special form):
$D = \Delta $, $v\left( {x,u} \right) = u$.

2) If $A$ is $u_{t} = a_{0} u^{\lambda} \left( {\Delta u + q_{0} u} \right) 
+ q_{1} u$ up to shift $u \to u - u_{0} $, $\lambda \ne 0$, $a_{0} ,q_{0} 
,q_{1} = \const$: 
$Df = \beta ^{\lambda - 1}\left(\Delta \left( \beta f \right)+q_0 \beta f \right)$ 
for some function $\beta :X \to \mathbb R$, 
$v\left( {x,u} \right) = \beta ^{ - 1}\left( {x} \right)u$, $Q = q_{1} v$.

3) If $A$ is $u_t=a_{0} e^{\lambda u}\left( \Delta u + q_0 \right) + q_{1} $, 
$\lambda \ne 0$, $a_{0} ,q_{0} ,q_{1} = \const$: 
\noindent
$Df = e^{\lambda \beta} \left( {\Delta f + \Delta \beta + q_{0}}  \right)$ 
for some function $\beta :X \to \mathbb R$, 

\noindent
$v\left( {x,u} \right) = u - \beta \left( {x} \right)$, $Q = q_{1} $.

We will call such morphisms "\textit{canonical}".
In the category $\mathcal {PE}$ the canonical representative in any class of morphisms 
is defined uniquely up to diffeomorphism of manifold $Y$, and in the 
category $\mathcal {APE}$ it is defined uniquely up to conformal 
diffeomorphism of $Y$.

Further we restrict ourselves by the investigation of the canonical maps for 
the first variant, that is will look for such maps $p$ from the given 
Riemann manifold $X$ onto arbitrary Riemann manifolds $Y$, for which 
Laplacian on $X$ is projected to some operator on $Y$ (note that this 
canonical maps will be canonical for given $X$ in the cases (2), (3) too).

Note that isomorphic autonomous equations $B$, factorized given $A$, 
are distinguished only by
arbitrary transformations $v \to v'(y,v)$ and has the same projection $p:x
\to y(x)$
up to conformal diffeomorphism of $Y$. Therefore to find such projection
$p:X \to Y$ for canonical morphism is to find all autonomous morphisms from
this equivalence class.

\section{Factorizing of heat equation in ${\mathbb R}^3$}

Let $\mathcal {DAPE}$ be full subcategory of $\mathcal {APE}$, whose objects 
are autonomous parabolic equations of divergent shape: 
\[
u_t = c(x,u)^{-1}\diverg \left( k(x,u)\nabla u \right) +q(x,u),
\]
and morphisms are autonomous morphisms of the category $\mathcal {APE}$.

\textbf{Theorem 6.} 
Let $X$ is a connected region of ${\mathbb R}^3$ with Euclidean metric, 
$Y$ is a manifold without boundary, $A$ does not have form (2 or 3) 
from Theorem 5. 
Then $p$ define canonical morphism of $A$ in $\mathcal {DAPE}$ 
iff $p$ is restriction on $X$ of factorization
${\mathbb R}^3$ under some (may be discontinuous) group $G$ of isometries.

\section{Factorizing with dimension decrease by 1}

\textbf{Theorem 7.} 
Let $A$ does not have form (2 or 3) from Theorem 5, and 
(a) $p:X \to Y$ is a fibering; (b) $X$ and $Y$ are oriented; 
(c) $X$ is an open domain in complete Riemann space $\tilde {X}$; 
(d) $\dim Y = \dim X - 1$. 
Then $p$ define canonical morphism to $\mathcal {DAPE}$ iff the following 
conditions fulfilled:

\noindent
a) $p$ is a superposition of maps $p_{1} :X \to {Y}'$ and $p_{0} :{Y}' \to Y$; 

\noindent
b) $p_{1} :X \to {Y}'$ is a restriction on $X$ of the projection $\tilde {X} 
\to {{\tilde {X}} \mathord{\left/ {\vphantom {{\tilde {X}} {G_{1}} }} 
\right. \kern-\nulldelimiterspace} {G_{1}} }$, where $G_{1} $ is some 
1-parameter subgroup of group $\Isom \left( {\tilde {X}} \right)$
of all isometries of $\tilde {X}$;

\noindent
c) $p_{0} :{Y}' \to \tilde {Y}$ is isomeric covering (for the metric on 
${Y}'$, inherited from $X$);

\noindent
d) For the vector field $\eta $ generating group $G_{1} $, the function 
$\vartheta = \left\langle {\eta ,\eta}  \right\rangle $, defined on ${Y}'$, 
is projectible on $Y$.

\section{Factorizing with dimension decrease by 1: comparison with group analysis}

As it was shown in the section 5, when we factorizing heat equation in 
${\mathbb R}^3$ with Euclidean metric, then the class of modeled 
(3D) solutions of A coincides with a class of solutions of A, which are 
invariant under some (maybe discontinuous) group of isometries of ${\mathbb R}^3$.

But this results about coincidence of factorizing maps 
for the heat equation in ${\mathbb R}^3$ with Euclidean metric with factormaps 
by symmetry groups (that is isometries groups) are accidental. 

At first, projection $p_{0} :{Y}' \to Y$ from previous section is not necessarily generated by some group of
transformation of ${Y}'$. 

At second, let even ${Y}' = {{Y} \mathord{\left/ {\vphantom {{Y} 
{G_{0}} }} \right. \kern-\nulldelimiterspace} {G_{0}} }$, 
where $G_{0} $ is some discrete group of the isometries of ${Y}'$.
The question is: could group $G_{0}$ be lifted to some group 
of the isometries of $X$, which preserve projection onto $Y$?

Let the group  $G_{1}$ be fixed, which satisfies conditions 
of Theorem 7. 
We consider differential-geometric connection $\chi $  on a fibering  
$p_{1} :X \to {Y}'$ with the structural group  $G_{1} $, which horizontal 
planes are orthogonal to  $G_{1} $ orbits.

\textbf{Theorem 8 (necessary condition).} 
If a discrete group  $G_{0} $, which 
operates on  ${Y}'$ and satisfies conditions of the 
Theorem 7, 
could be lifted to the subgroup of $\Isom (X)$, 
then curvature form $d\chi $, projected on  ${Y}'$,
would be invariant respectively $G_{0} $.

\textbf{Lemma 1.}
$\chi $ may be decomposed on a sum $\chi = p_{1 *}  {\chi} ' 
+ dh$, where ${\chi} ' \in T^{\ast} {Y}'$, $h$ is a function from $X$ to 
$\mathcal H$, $\mathcal H$ is fiber of $p_{1} $ (that is either $\mathbb R$, or circle 
$\mathbb R \bmod H$, where 
$H = \const$ is integral $\chi $ on a vertical cycle).

\textbf{Theorem 9 (necessary and sufficient condition).}
A discrete group $G_{0} $, operating on ${Y}'$ and satisfying conditions 
of the Theorem 7, could be lifted to the subgroup 
of $\Isom \left( {X} \right)$, iff
$\forall g \in G_{0}$ the form $g \chi ' - \chi '$ is:

- Exact, if the fiber of $p_{1} $ is simply connected;

- Closed with periods, multiply $H$, if the fiber of $p_{1} $ is multiply 
connected.

Particularly, if $X = \mathbb R^{n}$, and $G_{1} $ is the rotations group, $\eta = 
\sum\nolimits_{i = 1}^{m} {a_{i} \partial _{\varphi _{i}} }  $, $m \ge 3$, 
or $G_{1} $ is the screw motions group, $\eta = \partial _{z} + 
\sum\nolimits_{i = 1}^{m} {a_{i} \partial _{\varphi _{i}} }  $, $m \ge 2$, 
then such groups $G_{0} $ exist, which does not lift on $X$.

\section{Factorizing with dimension decrease}

Let's equip $X$ with connection generated by planes orthogonal to fibers.

\textbf{Theorem 10.} 
Let (a) $p:X \to Y$ is a fibering; (b) $\dim Y < \dim X$. 
Then $p$ define canonical morphism to $\mathcal {DAPE}$ 
iff the following conditions fulfilled:

\noindent
1) The fibers of $p$ are parallel;

\noindent
2) The transformation of a fiber over an initial point to a fiber over a 
final point changes volumes proportionally when we translate along any curve 
on $Y$;

\noindent
3) The holonomy group saves volume on a fiber.

\noindent
Moreover, $p$ define canonical morphism to $\mathcal {APE}$ 
iff conditions (1-2) fulfilled.

\textbf{Example 1}
($\dim X = 4$, $\dim Y = 2$).

\noindent
Let $X = \left\{ {\left( {x,y,z,w} \right)} \right\}$ with the metric
\[
g_{ij} = \left( {{\begin{array}{*{20}c}
 {1} \hfill & {0} \hfill & {0} \hfill & {0} \hfill \\
 {0} \hfill & {1 + \alpha ^{2} + \beta ^{2}} \hfill & {\alpha}  \hfill & 
{\beta}  \hfill \\
 {0} \hfill & {\alpha}  \hfill & {1} \hfill & {0} \hfill \\
 {0} \hfill & {\beta}  \hfill & {0} \hfill & {1} \hfill \\
\end{array}} } \right),
\alpha = xe^{w}, \quad \beta = xe^{z},
\]
\noindent
$Y = \left\{ {\left( {x,y} \right)} \right\}$ with the Euclidean metric, 
$p\left( {x,y,z,w} \right) = \left( {x,y} \right)$. 
Then  map $p$ and equation $v_t=v_{xx}+v_{yy}$ are factorization of 
the equation
\[
u_{t} = u_{xx} + u_{yy} - 2\alpha u_{yz} - 
\]
\[
-2\beta u_{yw} + \left( {1 + 
\alpha ^{2}} \right)u_{zz} +2\alpha \beta u_{zw} + \left( {1 + \beta ^{2}} 
\right)u_{ww} + 
\left( {\alpha \beta}  \right)_{w} u_{z} + \left( {\alpha 
\beta}  \right)_{z} u_{w} ,
\]
\noindent
where $\alpha = xe^{w}$ and $\beta = xe^{z}$, by the map $p:\left( {x,y,z,w} 
\right) \to \left( {x,y} \right)$. (The same is true for he equations 
$v_t=a(v)\Delta v$ on $Y$ and $u_t=a(u)\Delta u$ on $X$ for arbitrary 
function $a$, but for simplicity we will write linear equations in examples.)
However the only transformations $X$, 
under which both the last equation and all it's solutions 
projected by $p$ are invariant, are 
$\left( {x,y,z,w} \right) \to \left( {x,y,w,z} \right)$ and 
identity. Moreover, another transformation with such properties does not 
exist even locally (i.e. it couldn't be defined in any small neighborhood on 
$X$), even if we replace the request ``to keep the equation invariant'' by 
the request ``to be conformal''.

\textbf{Example 2}
($\dim X = 3$, $\dim Y = 2$).

\noindent
Let $\tilde {X} = 
\mathbb R^{3} = \left\{ {\left( {x,y,z} \right)} \right\}$ with the metric 
\[
g_{ij} = \left( {{\begin{array}{*{20}c}
 {1 + z^{2}} \hfill & {z} \hfill & { - z} \hfill \\
 {z} \hfill & {2} \hfill & { - 1} \hfill \\
 { - z} \hfill & { - 1} \hfill & {1} \hfill \\
\end{array}} } \right),
\]
\noindent
$\tilde {Y} = \left\{ {\left( {x,y} \right)} \right\}$ with the Euclidean 
metric. Let's consider group $H$ of isometries $\tilde {X}$, generated by 
the screw motion $\left( {x,y,z} \right) \to \left( {x + 1, - y, - z} 
\right)$ ($H$ is projectible on $\tilde {Y}$), $X = {{\tilde {X}} 
\mathord{\left/ {\vphantom {{\tilde {X}} {H}}} \right. 
\kern-\nulldelimiterspace} {H}}$, $Y = {{\tilde {Y}} \mathord{\left/ 
{\vphantom {{\tilde {Y}} {H}}} \right. \kern-\nulldelimiterspace} {H}}$, 
$p\left( {x,y,z} \right) = \left( {x,y} \right)$. $Y$ is homeomorphic to the 
Mobius band without a boundary; $X$ is homeomorphic to the torus without 
a boundary.

Then map $p$ and equation 
$e^{x}v_{t} = \left( {e^{x}v_{x}}  \right)_{x} + \left({e^{x}v_{y}}  \right)_{y} $, 
or $v_{t} = v_{xx} + v_{yy} + v_{x} $
on $Y$ are factorizations of the equation 
\[
u_{t} = u_{xx} + u_{yy} + u_{x} + 2zu_{xz} 
+ 2u_{yz} + \left( {\left( {2 + z^{2}} \right)u_{z}}  \right)_{z}
\]
\noindent
on $X$. 
However the only transformation $X$, under which both the last equation and all 
projected by $p$ it's solutions are invariant, is identity map. Moreover, 
there doesn't exist a non-identity conformal transformation $X$, under which 
all projected by $p$ solutions of the last equation are invariant.

\textbf{Example 3}
 ($\dim X = 3$, $\dim Y = 1$).

\noindent
Let $X = S^{1} \times \mathbb R^{2} = \left\{ {\left( {x,y,z} \right):x 
\in \mathbb R \bmod 1,y,z \in \mathbb R} \right\}$, 
equipped with the metric
\[
g_{ij} = \left( {{\begin{array}{*{20}c}
 {\alpha ^{2} + \beta ^{2}} \hfill & {\alpha}  \hfill & {\beta}  \hfill \\
 {\alpha}  \hfill & {1} \hfill & {0} \hfill \\
 {\beta}  \hfill & {0} \hfill & {1} \hfill \\
\end{array}} } \right),
\quad
\alpha = - e^{z},
\quad
\beta = 2y,
\]
\noindent
$\tilde {Y} = S^{1} = \left\{ {x \in \mathbb R \bmod 1} \right\}$ equipped
with the Euclidean metric, $p\left( {x,y,z} \right) = x$.
Then map $p$ and equation $v_{t} = v_{xx} $ on $Y$ are factorizations 
of the equation 
\[
u_{t} = u_{xx} + 
\]
\[
+\left( {1 + \alpha ^{2}} \right)u_{yy} + \left( {1 + 
\beta ^{2}} \right)u_{zz} + 
2\alpha \beta u_{yz} - 2\alpha u_{xy} - 2\beta u_{xz} + \left( {\alpha \beta}
\right)_{y} u_{z} + \left( {\alpha \beta}  \right)_{z} u_{y} 
\]
\noindent
on $X$. However the only transformation $X$, under which both 
the last equation and all 
projected by $p$ it's solutions are invariant, is identity map.

\textbf{Example 4}
($\dim X = 2$, $\dim Y = 1$). 

\noindent
Let $X = {{\mathbb R^{2}}
\mathord{\left/ {\vphantom {{\mathbb R^{2}} {G}}} \right. 
\kern-\nulldelimiterspace} {G}}$ 
with the Euclidean metric, when $G$ is the group generated by the 
sliding symmetry respectively the straight line $l$. The orthogonal 
projection of $X$ onto the mean circumference (image of the line $l$) define 
equation $v_{t} = v_{yy} $ on $l$, 
factorized the equation $u_{t} = u_{xx} + u_{yy}$
on $X$. However the only transformation $X$, under which both the last equation 
and all projected by $p$ it's solutions are invariant, is reflection 
respectively $l$.

\section{Factorization without dimension decrease}

If $\dim X = \dim Y$, then $p:X \to \tilde {Y}$ projected Laplacian 
iff it is isometric projection up to some conformal transformation $Y$. 

\textbf{Example 5.}
Let's manifold $X$ be a plane without 3 points: 
$A\left( {0,0} \right)$, $B\left( {1,0} \right)$ and $C\left( {0,2} \right)$. 
Let's consider heat equation on $X$ with metric 
$g_{ij} = \lambda ^{2}\left( {x} \right)\delta _{ij} $:
\begin{equation}
\label{Prokhorova:eq_Ex1}
\lambda ^{2}\left( {x} \right)u_{t} = u_{11} + u_{22} ,
\end{equation}
\noindent
where $\lambda \left( {x} \right) = \rho \left( {x,A} \right)\rho \left( 
{x,B} \right)\rho \left( {x,C} \right)$, $\rho $ is the distance function 
(in usual plane metric). 
Let $Y = X$, and map $p:X \to Y$ is given by the formula 
$y = \frac{1}{4} x^{4} - \frac{{1 + 2i}}{{3}}x^{3} + ix^{2}$, 
where $x$, $y$ are considered as points at a complex plane.

Because of $\left| {y_{x}}  \right| = \left| {x\left( {x - 1} \right)\left( 
{x - 2i} \right)} \right| = \lambda \left( {x} \right)$, heat equation 
$u_{t} = u_{11} + u_{22} $ on $Y$, equipped by Euclidean metric 
$g_{ij} = \delta _{ij} $, is factorisation of the equation (4) on the manifold 
$\mathop {X}\limits^{\circ } $, which is obtained by deleting 
of pre-images of images of zeroes of $\lambda $ from $X$. 
However, there does not exist non-identical transformation 
of $\mathop {X}\limits^{\circ }$, 
under which all projected by $p$ solutions of equation (4) 
are invariant. Moreover, there does not exist a non-identical transformation 
of any manifold $X'$, under which an equation (4) is 
invariant, if $X'$ is obtained by deleting an arbitrary discrete set of 
points from $X$.

\textbf{Example 6.}
Let's consider an equation on $X = \mathbb R^{2}$:
\begin{equation}
\label{Prokhorova:eq_Ex2}
u_{t} = \left( {1 + \left| {x} \right|^{2}} \right)^{2}\left( {u_{11} + 
u_{22}}  \right).
\end{equation}
\noindent
Let $g$ be the transformation of ${{\mathbb R^{2}} \mathord{\left/ {\vphantom 
{{\mathbb R^{2}} {\left\{ {0} \right\}}}} \right. \kern-\nulldelimiterspace} 
{\left\{ {0} \right\}}}$, that maps $x \in X$ to the point, obtained 
from $x$ by inversion under the unit circle with a center in an origin and 
consequent reflection under this center. Equation (5) 
is invariant with respect to $g$, but $g$ is not defined at origin. However the map 
$p:X \to Y = \mathbb P^{2}$ onto the projective plane, which past together points 
$x$ and $gx$ at $x \ne 0$, is defined on all $X$ and gives smooth 
projection. Then inducing on $Y$ heat equation is factorization of 
original equation on $X$.

\end{document}